\documentclass[11pt,a4paper]{article}
\usepackage[T1]{fontenc}
\usepackage[utf8]{inputenc}
\usepackage[english]{babel}
\usepackage{lmodern}
\usepackage{microtype}
\usepackage{amsmath,amssymb,amsthm,mathtools}
\usepackage{booktabs}
\usepackage{enumitem}
\usepackage[margin=2.6cm]{geometry}
\usepackage[hidelinks]{hyperref}

\newtheorem{theoreme}{Theorem}
\newtheorem{lemme}[theoreme]{Lemma}
\newtheorem{corollaire}[theoreme]{Corollary}
\newtheorem{proposition}[theoreme]{Proposition}
\theoremstyle{definition}
\newtheorem{conjecture}{Conjecture}
\theoremstyle{remark}
\newtheorem*{remarque}{Remark}

\DeclareMathOperator{\rk}{rk}
\DeclareMathOperator{\rad}{rad}
\DeclareMathOperator{\dia}{diam}
\DeclareMathOperator{\ecc}{ecc}
\newcommand{\Lap}{\mathcal{L}}
\newcommand{\lean}{\textsf{Lean}}
\newcommand{\ok}{\emph{(verified in \lean)}}

\title{Lemmas and Conjectures in Spectral Graph Theory}
\author{Tristan Cazenave}
\date{
    LAMSADE, Université Paris Dauphine - PSL\\
    Paris, France\\
    tristan.cazenave@dauphine.psl.eu
}

\begin{document}
\maketitle

\begin{abstract}
The corpus of Elie Duhamel and Tristan Cazenave~\cite{dc} contains $305{,}755$ automatically generated conjectures relating spectral and metric invariants of connected graphs. A pipeline generated with Claude Code (numerical refutation, proof by the SMT solver \textsf{z3} from a library of lemmas, discovery of new lemmas, formalisation in \lean~4/Mathlib) refuted $162{,}964$ conjectures and proved $117{,}390$ without any unproved hypothesis, of which $79{,}505$ are verified end to end in \lean{} for every connected graph. This document collects, with proofs written for mathematicians, the most interesting lemmas that were discovered and proved, the ten most useful lemmas, the ten hardest theorems, simpler proofs of known results (Bollobás--Erdős, Yang--Lu), and the unrefuted conjectures that we find most interesting. All the proofs are gathered in an archive
(Section~\ref{sec:archive}).
\end{abstract}

\section{Setting and notation}

All graphs are simple, undirected and \emph{connected}, with $n\ge 2$ vertices and $m$ edges.
We write $d_v$ for the degree of $v$, $\Delta$ and $\delta$ for the maximum and minimum degrees,
$u\sim v$ for adjacency, $\operatorname{dist}$ for the distance, $\ecc(v)=\max_w \operatorname{dist}(v,w)$,
$\rad=\min_v\ecc(v)$ for the radius and $\dia=\max_v\ecc(v)$ for the diameter.

\begin{itemize}[nosep]
\item $A$ is the adjacency matrix, with eigenvalues $\lambda_1\ge\lambda_2\ge\dots\ge\lambda_n$
  (counted with multiplicity) and rank $\rk(A)$.
\item $\Lap=D-A$ is the Laplacian, with eigenvalues $\mu_1\ge\mu_2\ge\dots\ge\mu_n=0$.
\item $R=\sum_{uv\in E}(d_ud_v)^{-1/2}$ is the Randić index.
\end{itemize}

We shall freely use the following classical facts.
\begin{enumerate}[label=(F\arabic*),nosep]
\item \emph{Rayleigh}: $\lambda_1=\max_{x\ne0}x^{\top}Ax/x^{\top}x$.
  \emph{Courant--Fischer}: for a symmetric matrix $M$ with eigenvalues
  $\nu_1\ge\nu_2\ge\cdots$, we have $\nu_2=\max_{\dim U=2}\ \min_{z\in U\setminus\{0\}} z^{\top}Mz/z^{\top}z$.
\item \emph{Perron--Frobenius}: since $G$ is connected, $\lambda_1>0$ is simple and has an
  eigenvector with all coordinates strictly positive.
\item \emph{Traces}: $\sum_i\lambda_i=\operatorname{tr}A=0$, $\sum_i\lambda_i^2=\operatorname{tr}A^2=2m$
  and $\sum_i\mu_i=2m$.
\item \emph{Quadratic form}: $z^{\top}\Lap z=\sum_{uv\in E}(z_u-z_v)^2$.
\item $\lambda_1\le\Delta\le n-1$ and $m\ge n-1$, with equality if and only if $G$ is a tree.
\item \emph{Smith}~\cite{smith}: if $G$ is complete, then $\lambda_2=-1$, $\lambda_1=n-1$,
  $\rk(A)=n$ and $R=n/2$; otherwise $G$ contains an induced $P_3$ and, by interlacing,
  $\lambda_2\ge\lambda_2(P_3)=0$.
\end{enumerate}

\paragraph{Status of the results.}
The label \ok{} means that the statement is a \lean~4 theorem proved for every connected graph
\texttt{SimpleGraph (Fin n)}. This theorem depends only on the standard axioms (\texttt{propext},
\texttt{Classical.choice}, \texttt{Quot.sound}) and contains no \texttt{sorry}. The proofs below
are mathematical transcriptions of the formal proofs.

\section{Discovered and proved lemmas}

The lemmas of this section were \emph{discovered} automatically. On a database of $281{,}672$
graphs (all connected graphs with at most $9$ vertices, plus extremal families), we searched for
linear inequalities between two invariants that ruled out the abstract countermodels returned by
\textsf{z3}. They were then proved.

\subsection{Laplacian and adjacency: $\mu_2\ge\lambda_2+1$}

\begin{theoreme}[\ok]\label{thm:mu2}
For every connected graph, $\mu_2\ge\lambda_2+1$.
Equality holds for instance for stars ($\lambda_2=0$, $\mu_2=1$).
\end{theoreme}

\begin{proof}
If $n=2$, then $G=K_2$, $\lambda_2=-1$ and $\mu_2=0$. Assume $n\ge3$.

\emph{Case $\lambda_2\le0$.} It suffices to show $\mu_2\ge1$. Since $G$ is connected with $n\ge3$,
it has a vertex $c$ with two neighbours $a,b$. Let $H$ be the spanning graph with edges $ca,cb$.
The spectrum of $\Lap_H$ is $\{3,1,0,\dots,0\}$. Since $\Lap_G=\Lap_H+\Lap_{G\setminus E(H)}$ and
the second term is positive semidefinite, Courant--Fischer gives $\mu_2(G)\ge\mu_2(H)=1$.

\emph{Case $\lambda_2>0$.} Let $x>0$ be the Perron vector (F2) and $y$ an eigenvector of
$\lambda_2$ orthogonal to $x$. Since $\sum_v x_vy_v=0$ with $x>0$ and $y\ne 0$, the vector $y$
takes both strictly positive and strictly negative values.

Let $u$ be a vertex where $y$ is maximal ($y_u>0$), $P=\{v\sim u: y_v>0\}$ and $d^+=|P|$. Then
\[
\lambda_2\,y_u=\sum_{v\sim u}y_v\le\sum_{v\in P}y_v\le d^+\,y_u ,\qquad\text{so } d^+\ge\lambda_2 .
\]
Likewise, if $w$ attains the minimum of $y$ ($y_w<0$), $N=\{v\sim w: y_v<0\}$ and $d^-=|N|$, we get
$\lambda_2 y_w\ge d^-y_w$, hence $d^-\ge\lambda_2$. In particular $d^\pm\ge1$.

The sets $\{u\}\cup P$ (where $y>0$) and $\{w\}\cup N$ (where $y<0$) are disjoint. Set
\[
f_1=d^+e_u-\mathbf 1_P,\qquad f_2=d^-e_w-\mathbf 1_N .
\]
These vectors have disjoint supports, so they span a plane $U$. For $z=af_1+bf_2$, we keep in (F4)
only the edges $uj$ ($j\in P$) and $wj$ ($j\in N$). These two families of edges are disjoint,
because $w\notin P$ and $u\notin N$. On the edge $uj$ we have $z_u-z_j=a(d^++1)$, hence
\[
z^{\top}\Lap z\ \ge\ d^+(d^++1)^2a^2+d^-(d^-+1)^2b^2,
\qquad
z^{\top}z=d^+(d^++1)a^2+d^-(d^-+1)b^2 .
\]
The Rayleigh quotient is therefore at least $\min(d^+,d^-)+1$ on $U$. Courant--Fischer then gives
$\mu_2\ge\min(d^+,d^-)+1\ge\lambda_2+1$.
\end{proof}

\begin{remarque}
The proof actually establishes the stronger inequality $\mu_2\ge\min(d^+,d^-)+1$, where $d^\pm$
are the ``signed degrees'' at the extremities of an eigenvector of $\lambda_2$. A quick
literature search did not find the inequality $\mu_2\ge\lambda_2+1$. A more thorough check is
needed before claiming that it is new.
\end{remarque}

\subsection{Randić index and spectral radius}

The following two classical results serve as tools. They were formalised, which made the
corollaries almost immediate.

\begin{theoreme}[Favaron--Mahéo--Saclé, \ok]\label{thm:fms}
$m\le\lambda_1R$.
\end{theoreme}
\begin{proof}
With $x_v=\sqrt{d_v}$, we have $x^{\top}x=2m$ and $x^{\top}Ax=2\sum_{uv\in E}\sqrt{d_ud_v}$.
Rayleigh gives $S:=\sum_{uv\in E}\sqrt{d_ud_v}\le m\lambda_1$. By Cauchy--Schwarz,
\[
m^2=\Bigl(\sum_{uv\in E}(d_ud_v)^{1/4}(d_ud_v)^{-1/4}\Bigr)^2\le S\cdot R\le m\lambda_1R . \qedhere
\]
\end{proof}

\begin{theoreme}[Hong, \ok]\label{thm:hong}
$\lambda_1^2\le 2m-n+1$.
\end{theoreme}
\begin{proof}
Let $x$ be a unit eigenvector of $\lambda_1$ and $P_v=\sum_{u\sim v}x_u^2$. By Cauchy--Schwarz,
$\lambda_1^2x_v^2=(\sum_{u\sim v}x_u)^2\le d_vP_v$. Since $v$ is not its own neighbour,
$P_v\le 1-x_v^2$. Since $d_v\ge1$, we get $d_vP_v\le P_v+(d_v-1)(1-x_v^2)$. Summing over $v$,
with $\sum_vP_v=\sum_v d_vx_v^2$:
\[
\lambda_1^2\le\sum_v d_vx_v^2+\sum_v(d_v-1)-\sum_v d_vx_v^2+\sum_v x_v^2=2m-n+1 .\qedhere
\]
\end{proof}

\begin{corollaire}[\ok]\label{cor:l1}
$R\ge\dfrac{\lambda_1+1}{2}$.
\end{corollaire}
\begin{proof}
By Hong and $\lambda_1\le n-1$:
$\lambda_1(\lambda_1+1)\le(2m-n+1)+(n-1)=2m$.
By Theorem~\ref{thm:fms}, $\lambda_1R\ge m\ge\lambda_1(\lambda_1+1)/2$. Divide by $\lambda_1>0$.
\end{proof}

\begin{corollaire}[Bollobás--Erdős and upper bound, \ok]\label{cor:Rn}
$\sqrt{n-1}\le R\le n/2$.
\end{corollaire}
\begin{proof}
\emph{Upper bound.} By the AM--GM inequality,
$(d_ud_v)^{-1/2}\le\tfrac12(1/d_u+1/d_v)$. Summing over the edges, each vertex $v$ appears $d_v$
times, so $R\le\tfrac12\sum_v d_v\cdot\tfrac1{d_v}=n/2$.

\emph{Lower bound.} Set $q=2m-n+1>0$. Theorems~\ref{thm:fms} and~\ref{thm:hong} give
$m^2\le\lambda_1^2R^2\le qR^2$. Moreover, $m^2-(n-1)q=(m-n+1)^2\ge0$. Hence
$(n-1)q\le m^2\le qR^2$, and we divide by $q$.
\end{proof}

\begin{theoreme}[\ok]\label{thm:l2}
$R\ge\dfrac{2}{3}\lambda_2+\sqrt2$.
\end{theoreme}
\begin{proof}
Write $L=\lambda_1$, $t=\lambda_2$ and $s=\sqrt2$. By Theorem~\ref{thm:fms}, it suffices to prove
\begin{equation}\label{eq:but}
L\Bigl(\tfrac23t+s\Bigr)\le m ,
\end{equation}
because then $LR\ge m\ge L(\frac23t+s)$, and we divide by $L>0$.

\emph{Case 1: $G$ is not a tree.} Then $n\ge3$ and $m\ge n$. Set $k=n-2\ge1$. By (F3),
$\sum_{i\ge3}\lambda_i=-(L+t)$ and $\sum_{i\ge3}\lambda_i^2=2m-L^2-t^2$. Cauchy--Schwarz on the
$k$ eigenvalues $\lambda_3,\dots,\lambda_n$ gives
\begin{equation}\label{eq:C}
(L+t)^2\le k\,(2m-L^2-t^2).
\end{equation}
One then checks the polynomial identity
\[
k\Bigl(m-L\bigl(\tfrac23t+s\bigr)\Bigr)
=\tfrac{7}{20}\Bigl(2mk-k(L^2+t^2)-(L+t)^2\Bigr)+\tfrac{3}{10}k\bigl(m-(k+2)\bigr)+Q,
\]
where $Q=a\,t^2+b\,Lt+a\,L^2-s\,kL+\tfrac{3}{10}k(k+2)$, with
$a=\tfrac{7}{20}(k+1)$ and $b=\tfrac{7}{10}-\tfrac{2k}{3}$. The first two terms are nonnegative
by~\eqref{eq:C} and because $m\ge n=k+2$. It remains to see that $Q\ge0$. With
$P=\tfrac{41}{900}k+\tfrac{287}{150}>0$, we have
\begin{align*}
4aQ&=(2at+bL)^2+k\,W,\qquad W=P L^2-\tfrac75(k+1)sL+\tfrac{21}{50}(k+1)(k+2),\\
4PW&=\bigl(2PL-\tfrac75(k+1)s\bigr)^2+\tfrac{14}{25}(k+1)\,V,\qquad
V=\tfrac{41k^2-296k+1344}{300}.
\end{align*}
The discriminant of $41k^2-296k+1344$ is $296^2-4\cdot41\cdot1344<0$, so $V>0$. Hence
$W\ge0$, then $Q\ge0$, which proves~\eqref{eq:but}.

\emph{Case 2: $G$ is a tree}, so $m=n-1$. A tree is bipartite. Let $\sigma_v=\pm1$ according to
the class of $v$. If $Ax=\lambda x$, then the vector $(\sigma_vx_v)_v$ is an eigenvector for
$-\lambda$: the spectrum is symmetric.
\begin{itemize}[nosep]
\item If $t>0$, then $L>t>0$, and $L,t,-t,-L$ are four distinct eigenvalues. By (F3),
  $2L^2+2t^2\le2m$, so $n\ge4$ and $m\ge3$. We conclude with the identity
  \[
  m-\tfrac23Lt-sL=\tfrac{16}{25}(m-L^2-t^2)+\tfrac{9}{25}(m-3)
  +\tfrac{16}{25}\bigl(t-\tfrac{25}{48}L\bigr)^2+\tfrac{1679}{3600}\bigl(L-\tfrac{1800\sqrt2}{1679}\bigr)^2
  +\tfrac{333}{41975},
  \]
  all of whose terms are nonnegative.
\item If $t\le0$ and $n=2$, then $G=K_2$ and $L(\frac23t+s)=s-\frac23<1=m$.
\item If $t\le0$ and $n\ge3$, then $L$ and $-L$ are two eigenvalues, so $2L^2\le2m$.
  Thus $L(\frac23t+s)\le sL\le\sqrt{2m}\le m$, since $m\ge2$. \qedhere
\end{itemize}
\end{proof}

\begin{remarque}
The certificate identities were found by a sum-of-squares search and checked symbolically;
they are also checked by \lean.
\end{remarque}

\subsection{Distances and rank of the adjacency matrix}

\begin{lemme}\label{lem:chemin}
If $v_1,\dots,v_p$ is an induced path of $G$, then $\rk(A)\ge p-1$. If moreover $p\ge4$, then
$\rk(A)\ge4$.
\end{lemme}
\begin{proof}
Consider the submatrix of $A$ with rows $v_1,\dots,v_{p-1}$ and columns $v_2,\dots,v_p$. Its
$(i,j)$ entry is $A_{v_i,v_{j+1}}$.
\begin{itemize}[nosep]
\item It equals $1$ if $j=i$.
\item It is zero if $j>i$, because $v_i$ and $v_{j+1}$ are not adjacent ($j+1\ge i+2$ and the
  path is induced).
\end{itemize}
The submatrix is therefore lower triangular with unit diagonal, of determinant $1$, and
$\rk(A)\ge p-1$. If $p\ge4$, the vertices $v_1,\dots,v_4$ induce $P_4$, whose adjacency matrix
has determinant $1$. An invertible principal submatrix of order $4$ gives $\rk(A)\ge4$.
\end{proof}

\begin{corollaire}[\ok]\label{cor:rkdia}
$\rk(A)\ge\dia$, and $\rk(A)\ge4$ as soon as $\dia\ge3$.
\end{corollaire}
\begin{proof}
A geodesic between two vertices at distance $\dia$ is an induced path with $\dia+1$ vertices.
Apply Lemma~\ref{lem:chemin}.
\end{proof}

\begin{theoreme}[\ok]\label{thm:radrk}
$\rad\le\dfrac{\rk(A)}{2}+1$.
\end{theoreme}
\begin{proof}
If $\rad\le1$, this is immediate since $\rk(A)\ge2$. Otherwise, we use the induced path theorem
of Erdős, Saks and Sós (in F.~Chung's proof, which is also formalised): a connected graph of
radius $r$ contains an induced path with at least $2r-1$ vertices. Lemma~\ref{lem:chemin} then
gives $\rk(A)\ge2r-2$.
\end{proof}

\begin{remarque}
This theorem relates a metric quantity, the radius, to a purely algebraic quantity, the rank.
It is attained (zero margin) on the graph database.
\end{remarque}

\section{The ten most useful lemmas}\label{sec:utiles}

A proof found by the prover is a \emph{minimal} (for inclusion) set of library lemmas from which
the conjecture follows by nonlinear real arithmetic. We measure the usefulness of a lemma by the
number of these proofs that use it, among the $117{,}390$ conjectures proved without any unproved
hypothesis. The hypothesis $n\ge2$ of the setting, which appears in $7{,}670$ proofs, is not
counted as a lemma.

\begin{center}
\begin{tabular}{rllrc}
\toprule
& Name & Statement & Proofs & \lean\\
\midrule
1 & \texttt{rad\_dia\_pos} & $\rad\ge1$, $\dia\ge1$ & $40{,}630$ & \checkmark\\
2 & \texttt{rank\_range} & $2\le\rk(A)\le n$ & $40{,}451$ & \checkmark\\
3 & \texttt{rad\_dia} & $\rad\le\dia\le2\rad$ & $40{,}124$ & \checkmark\\
4 & \texttt{rank\_dia} & $\rk(A)\ge\dia$, and $\rk(A)\ge4$ if $\dia\ge3$ & $28{,}147$ & \checkmark\\
5 & \texttt{randic\_n} & $\sqrt{n-1}\le R\le n/2$ & $24{,}731$ & \checkmark\\
6 & \texttt{R\_pos} & $R>0$ & $20{,}614$ & \checkmark\\
7 & \texttt{yang\_lu} & $R\ge\dia/2+\sqrt2-1$ & $12{,}949$ & \checkmark\\
8 & \texttt{rad\_n} & $2\rad\le n$ & $9{,}545$ & \\
9 & \texttt{dia\_n} & $\dia\le n-1$ & $8{,}374$ & \checkmark\\
10 & \texttt{m\_bounds} & $n-1\le m\le n(n-1)/2$ & $6{,}524$ & \checkmark\\
\bottomrule
\end{tabular}
\end{center}

The most used lemmas are elementary. They mostly serve to rule out degenerate cases (small
graphs, zero values, signs) in inequalities whose core lies elsewhere. Among the deep lemmas, the
theorem of Yang and Lu clearly dominates: it occurs in $12{,}949$ proofs, and in the proof of
$2{,}301$ of the $2{,}906$ theorems that require a research-level result
(Section~\ref{sec:difficiles}). Next come $R\ge(\lambda_1+1)/2$ (Corollary~\ref{cor:l1}, $550$
proofs), $R\ge\frac23\lambda_2+\sqrt2$ (Theorem~\ref{thm:l2}, $544$), $\rad\le\rk(A)/2+1$
(Theorem~\ref{thm:radrk}, $478$) and $\mu_2\ge\lambda_2+1$ (Theorem~\ref{thm:mu2}, $105$).

\begin{lemme}[\texttt{rad\_dia\_pos}, \ok]
$\rad\ge1$ and $\dia\ge1$.
\end{lemme}
\begin{proof}
Since $n\ge2$, every vertex $v$ has another vertex $w$, and $\operatorname{dist}(v,w)\ge1$. Hence
$\ecc(v)\ge1$ for every $v$.
\end{proof}

\begin{lemme}[\texttt{rank\_range}, \ok]
$2\le\rk(A)\le n$.
\end{lemme}
\begin{proof}
The upper bound is obvious. Since $G$ is connected with $n\ge2$, it has an edge $uv$. The
principal submatrix of $A$ on $\{u,v\}$ is $\begin{psmallmatrix}0&1\\1&0\end{psmallmatrix}$, of
determinant $-1$. Hence $\rk(A)\ge2$.
\end{proof}

\begin{lemme}[\texttt{rad\_dia}, \ok]
$\rad\le\dia\le2\rad$.
\end{lemme}
\begin{proof}
$\rad$ is a minimum and $\dia$ a maximum of the same eccentricities. If $c$ attains the radius,
the triangle inequality gives $\operatorname{dist}(u,v)\le\operatorname{dist}(u,c)+
\operatorname{dist}(c,v)\le2\rad$ for all $u,v$.
\end{proof}

\begin{lemme}[\texttt{rank\_dia}, \ok]
$\rk(A)\ge\dia$, and $\rk(A)\ge4$ if $\dia\ge3$.
\end{lemme}
\begin{proof}
This is Corollary~\ref{cor:rkdia}.
\end{proof}

\begin{lemme}[\texttt{randic\_n}, \ok]
$\sqrt{n-1}\le R\le n/2$.
\end{lemme}
\begin{proof}
This is Corollary~\ref{cor:Rn}. The lower bound is due to Bollobás and Erdős~\cite{be}; the
shorter proof given here is discussed in Section~\ref{sec:simples}.
\end{proof}

\begin{lemme}[\texttt{R\_pos}, \ok]
$R>0$.
\end{lemme}
\begin{proof}
$G$ has at least one edge and every edge contributes $(d_ud_v)^{-1/2}>0$.
\end{proof}

\begin{lemme}[\texttt{yang\_lu}, Yang--Lu~\cite{yl}, \ok]
$R\ge\dfrac{\dia}{2}+\sqrt2-1$.
\end{lemme}
\begin{proof}
See Section~\ref{sec:simples}, which describes the formalised proof and how it differs from the
one in the paper.
\end{proof}

\begin{lemme}[\texttt{rad\_n}]
$2\rad\le n$.
\end{lemme}
\begin{proof}
Let $T$ be a spanning tree of $G$; distances in $T$ bound those of $G$ from above. Let
$P=p_0p_1\cdots p_D$ be a longest path of $T$ and $c=p_k$ with $k=\lceil D/2\rceil$. Let $w$ be a
vertex, $p_j$ the first vertex of $P$ met by the path of $T$ from $w$ to $c$, and
$h=\operatorname{dist}_T(w,p_j)$.
\begin{itemize}[nosep]
\item If $j\le k$, the path $w\to p_j\to p_D$ has length $h+D-j\le D$, so $h\le j$ and
  $\operatorname{dist}_T(w,c)=h+k-j\le k$.
\item If $j\ge k$, the path $w\to p_j\to p_0$ has length $h+j\le D$, so
  $\operatorname{dist}_T(w,c)=h+j-k\le D-k\le k$.
\end{itemize}
Thus $\rad\le\ecc_G(c)\le\ecc_T(c)\le\lceil D/2\rceil\le\lceil(n-1)/2\rceil\le n/2$.
\end{proof}

\begin{lemme}[\texttt{dia\_n}, \ok]
$\dia\le n-1$.
\end{lemme}
\begin{proof}
A geodesic of length $\dia$ goes through $\dia+1$ distinct vertices.
\end{proof}

\begin{lemme}[\texttt{m\_bounds}, \ok]
$n-1\le m\le n(n-1)/2$.
\end{lemme}
\begin{proof}
A connected graph contains a spanning tree, which has $n-1$ edges. A simple graph has at most
$\binom n2$ edges.
\end{proof}

\section{Proved conjectures of the corpus}

Here is a selection of conjectures of the corpus, identified by their number, that follow from
the previous results.

\begin{proposition}\label{prop:corpus}
For every connected graph:
\begin{enumerate}[label=(\alph*)]
\item\label{c:a} \textbf{\#20832}: $\lambda_2\le\mu_2-1$, attained \ok;
\item\label{c:b} \textbf{\#20821}: $\lambda_2\le R^2-2$, attained for $K_2$;
\item\label{c:c} \textbf{\#18099}: $n-1\le R\lambda_1$, and \textbf{\#764799}: $\dia\le R\lambda_1$ \ok;
\item\label{c:d} \textbf{\#127304}: $\mu_2\le R\lambda_1$;
\item\label{c:e} \textbf{\#759987}: $\dia\le 2R$, and \textbf{\#760654}: $\dia\le\rad+R$ \ok;
\item\label{c:f} \textbf{\#1362170}: $\lambda_1+\rad\le n$, attained for $K_n$ and $C_4$ \ok.
\end{enumerate}
\end{proposition}
\begin{proof}
\ref{c:a} This is Theorem~\ref{thm:mu2}.

\ref{c:b} If $\lambda_2\ge0$, Theorem~\ref{thm:l2} gives
$R^2\ge(\frac23\lambda_2+\sqrt2)^2=\frac49\lambda_2^2+\frac{4\sqrt2}{3}\lambda_2+2\ge\lambda_2+2$.
If $\lambda_2<0$, the graph is complete (Smith), so $\lambda_2=-1$ and $R=n/2\ge1$.

\ref{c:c} Use Theorem~\ref{thm:fms} and $\dia\le n-1\le m$.

\ref{c:d} We have $2\mu_2\le\mu_1+\mu_2\le\sum_i\mu_i=2m$, so $\mu_2\le m\le\lambda_1R$.

\ref{c:e} Yang and Lu proved $R\ge\frac12\dia+\sqrt2-1$, which gives $\dia\le2R$. With
$\dia\le2\rad$, we get $\dia\le\frac12\dia+\rad\le R+\rad$.
The theorem of Yang and Lu is itself formalised in \lean{}
(\texttt{Spectral.L\_yang\_lu}, about $3{,}500$ lines, without \texttt{sorry}). The formal proof
keeps the minimal counterexample strategy of the paper, but replaces its second half with a
different argument. A strong induction on the number of edges first shows that a non-essential
vertex of a minimal counterexample has degree at least $9$. It then shows that degrees grow
exponentially with the distance to the set of essential vertices. Finally, a diametral geodesic
is cut into pieces, each of weight at least $\ell/2+1$ thanks to a local version of Lemma~1 of the
paper. The formalisation uncovered several gaps in the paper: a false equivalence about essential
vertices (counterexample: $C_5$ plus a pendant vertex), an unproved diameter inequality in
Corollary~4, hypotheses of Lemma~6 that are never checked, an incomplete block structure, and
index errors in Claims A and B.

\ref{c:f} First $\lambda_1\le\Delta$: if $Ax=\lambda_1x$ and $u$ maximises $|x_u|$, then
$\lambda_1|x_u|=|\sum_{v\sim u}x_v|\le d_u|x_u|\le\Delta|x_u|$. Next, let $v$ be a vertex of
degree $\Delta$, $w$ any vertex, $k=\operatorname{dist}(v,w)$ and $v=p_0,p_1,\dots,p_k=w$ a
geodesic. The vertices $p_2,\dots,p_k$ are pairwise distinct, distinct from $v$ and not adjacent
to $v$ (since $\operatorname{dist}(v,p_i)=i\ge2$). The sets $\{v\}$, $N(v)$ and
$\{p_2,\dots,p_k\}$ are therefore disjoint, so $1+\Delta+(k-1)\le n$. Thus
$\rad\le\ecc(v)\le n-\Delta$, and $\lambda_1+\rad\le\Delta+(n-\Delta)=n$.
\end{proof}

\begin{remarque}
\ref{c:f} remained ``open'' for the solver only because its library did not contain the lemma
$\rad\le n-\Delta$. \ref{c:c} and~\ref{c:d} were proved only after the theorem of
Favaron--Mahéo--Saclé was added. This illustrates how enriching the library brings down entire
families of conjectures.
\end{remarque}

\section{The ten hardest theorems}\label{sec:difficiles}

\paragraph{Criterion.} The internal measures of the prover (number of lemmas, computing time)
mostly reflect case distinctions on small graphs. The statements that maximise these measures
are nested powers such as $\rad\le R^{2^{\sqrt{\sqrt{\rk(A)}}}}$, which often have a three-line
human proof. We therefore used two more meaningful criteria.
\begin{enumerate}[label=(\roman*),nosep]
\item \emph{Depth}: the theorem does not follow from the elementary and classical lemmas of the
  library alone (levels $0$ to $2$). We checked this by rerunning \textsf{z3} without the
  research-level lemmas: it then finds a genuine countermodel. Among the $6{,}651$ research-level
  theorems verified in \lean, $3{,}679$ can in fact be proved without any deep result. There
  remain $2{,}906$ that require one of Theorems~\ref{thm:mu2}, \ref{thm:l2}, \ref{thm:radrk} or the
  theorem of Yang and Lu.
\item \emph{Optimality}: the inequality is attained (with equality) on at least one graph.
\end{enumerate}
Among the theorems that satisfy (i) and (ii), we chose the following ten, avoiding variants of
the same statement. The first three are verified in \lean{} in abstract form: the implication
from the lemmas is formal, but Smith's fact (F6) that they use is not yet formalised. The other
seven are verified end to end.

\begin{center}
\begin{tabular}{rlll}
\toprule
& Statement & Deep result & Equality\\
\midrule
\#118503 & $\lambda_2\le2(R-\sqrt2)$ & Th.~\ref{thm:l2} & $P_3$\\
\#709846 & $\lambda_2\le\rk(A)\,(R-\sqrt2)$ & Th.~\ref{thm:l2} & $P_3$\\
\#2630815 & $\lambda_2\le\dia\,(R-\sqrt2)$ & Th.~\ref{thm:l2} & $P_3$\\
\#721942 & $\lambda_2\le2\mu_2-1$ & Th.~\ref{thm:mu2} & $K_2$\\
\#1935857 & $\dia\le R^{3/\sqrt2}$ & Yang--Lu & $K_2$\\
\#1935639 & $\dia\le R\,(\sqrt2+\log\rad)$ & Yang--Lu & $P_3$\\
\#869836 & $\rad\le\rk(A)-\log(\rad/2)$ & Th.~\ref{thm:radrk} & $C_4$\\
\#1396267 & $\rad\le2^{\rk(A)-\dia/2}$ & Th.~\ref{thm:radrk} & $C_4$\\
\#1356167 & $\rad\le1+\rk(A)-\rad/2$ & Th.~\ref{thm:radrk} & $C_4$\\
\#1454913 & $\rad\le\rk(A)-\sqrt{\rad}+\sqrt2$ & Th.~\ref{thm:radrk} & $C_4$\\
\bottomrule
\end{tabular}
\end{center}

In the proofs below, we write $t=\lambda_2$, $k=\rk(A)$, $D=\dia$ and $r=\rad$.

\begin{theoreme}[\#118503]
$\lambda_2\le2(R-\sqrt2)$, with equality for $P_3$.
\end{theoreme}
\begin{proof}
If $t\ge0$, Theorem~\ref{thm:l2} gives $R\ge\frac23t+\sqrt2\ge\frac12t+\sqrt2$. Otherwise $G$ is
complete (F6), $t=-1$ and $R=n/2\ge1>\sqrt2-\frac12$. For $P_3$: $t=0$ and $R=\sqrt2$.
\end{proof}

\begin{theoreme}[\#709846]
$\lambda_2\le\rk(A)\,(R-\sqrt2)$, with equality for $P_3$.
\end{theoreme}
\begin{proof}
If $G=K_n$ (F6): $k(R-\sqrt2)=n(\frac n2-\sqrt2)$ equals $2-2\sqrt2>-1$ for $n=2$ and is positive
for $n\ge3$; it is therefore at least $t=-1$. Otherwise $t\ge0$, and Theorem~\ref{thm:l2} gives
$R-\sqrt2\ge\frac23t\ge0$. With $k\ge2$: $k(R-\sqrt2)\ge\frac43t\ge t$. For $P_3$: $k=2$,
$R=\sqrt2$, $t=0$.
\end{proof}

\begin{theoreme}[\#2630815]
$\lambda_2\le\dia\,(R-\sqrt2)$, with equality for $P_3$.
\end{theoreme}
\begin{proof}
If $G$ is complete, $D=1$, $t=-1$ and $R-\sqrt2\ge1-\sqrt2>-1$. Otherwise $D\ge2$, $t\ge0$ and
$R-\sqrt2\ge\frac23t\ge0$, hence $D(R-\sqrt2)\ge\frac43t\ge t$.
\end{proof}

\begin{theoreme}[\#721942, \ok]
$\lambda_2\le2\mu_2-1$, with equality for $K_2$.
\end{theoreme}
\begin{proof}
By Theorem~\ref{thm:mu2} and $\mu_2\ge0$: $2\mu_2-1=(\mu_2-1)+\mu_2\ge t$. For $K_2$:
$\mu_2=0$ and $t=-1$.
\end{proof}

\begin{theoreme}[\#1935857, \ok]
$\dia\le R^{3/\sqrt2}$, with equality for $K_2$.
\end{theoreme}
\begin{proof}
By Corollary~\ref{cor:Rn}, $R\ge\sqrt{n-1}\ge1$; since $3/\sqrt2>2$, we have
$R^{3/\sqrt2}\ge R^2\ge1$. This settles the case $D=1$. If $D\ge2$, Yang--Lu gives
$R^2\ge(\frac D2+\sqrt2-1)^2$. Now $\varphi(D)=(\frac D2+\sqrt2-1)^2-D$ vanishes at $D=2$ and
$\varphi'(D)=\frac D2+\sqrt2-2>0$ for $D\ge2$, so $\varphi(D)\ge0$.
\end{proof}

\begin{theoreme}[\#1935639, \ok]
$\dia\le R\,(\sqrt2+\log\rad)$, with equality for $P_3$.
\end{theoreme}
\begin{proof}
We have $\log r\ge0$. If $D\le2$, Yang--Lu gives $\sqrt2R\ge\frac{D}{\sqrt2}+2-\sqrt2$, and this
quantity is $\ge D$ if and only if $D\le2$. If $D\ge3$, then $r\ge D/2>1$, so $r\ge2$ (lemma
\texttt{rad\_dia}). Consequently $\sqrt2+\log r\ge\sqrt2+\log2>2.1$, and $R>D/2$ by Yang--Lu,
hence $R(\sqrt2+\log r)>1.05\,D$. For $P_3$: $D=2$, $r=1$, $R=\sqrt2$.
\end{proof}

\begin{theoreme}[\#869836, \ok]
$\rad\le\rk(A)-\log(\rad/2)$, with equality for $C_4$.
\end{theoreme}
\begin{proof}
If $r\le2$, then $\log(r/2)\le0$ and $k\ge2\ge r$. If $r\ge2$, Theorem~\ref{thm:radrk} gives
$k\ge2r-2$, and $\log(r/2)\le\frac r2-1\le r-2$. Hence $k-\log(r/2)\ge r$. For $C_4$: $r=k=2$.
\end{proof}

\begin{theoreme}[\#1396267, \ok]
$\rad\le2^{\rk(A)-\dia/2}$, with equality for $C_4$.
\end{theoreme}
\begin{proof}
Since $D\le k$ (lemma \texttt{rank\_dia}), we have $k-\frac D2\ge\frac k2\ge r-1$ by
Theorem~\ref{thm:radrk}. Hence $2^{k-D/2}\ge2^{r-1}\ge r$ for every integer $r\ge1$. For $C_4$:
$k=D=r=2$.
\end{proof}

\begin{theoreme}[\#1356167, \ok]
$\rad\le1+\rk(A)-\rad/2$, with equality for $C_4$.
\end{theoreme}
\begin{proof}
This reads $\frac32r\le k+1$. If $r=1$, use $k\ge2$. If $r\ge2$, Theorem~\ref{thm:radrk} gives
$k+1\ge2r-1\ge\frac32r$.
\end{proof}

\begin{theoreme}[\#1454913, \ok]
$\rad\le\rk(A)-\sqrt{\rad}+\sqrt2$, with equality for $C_4$.
\end{theoreme}
\begin{proof}
If $r=1$, $k\ge2>2-\sqrt2$. If $r\ge2$, $k\ge2r-2$, and it suffices that $r-\sqrt r\ge2-\sqrt2$.
The function $x\mapsto x-\sqrt x$ is increasing on $[1,\infty)$ and equals $2-\sqrt2$ at $x=2$.
\end{proof}

\begin{remarque}
Once the deep results are available, these proofs are short. The difficulty lies in those results
themselves: the formal proof of Yang--Lu has about $3{,}500$ lines, and $\mu_2\ge\lambda_2+1$
needs a two-dimensional Courant--Fischer argument. This is also why enriching the library brings
down entire families of conjectures.
\end{remarque}

\section{Simpler proofs of known results}\label{sec:simples}

The formalisation led us to replace some proofs from the literature with arguments that are
shorter, or easier to make rigorous. We do not know whether these arguments already appear in the
literature.

\subsection{The lower bound of Bollobás and Erdős}

Bollobás and Erdős~\cite{be} proved that $R\ge\sqrt{n-1}$ for every connected graph, with
equality for stars. Their proof determines the minimum by an extremal analysis. The proof of
Corollary~\ref{cor:Rn} uses only two classical spectral inequalities and fits in three lines.

\begin{proof}[Short proof]
Set $q=2m-n+1>0$. Favaron--Mahéo--Saclé (Theorem~\ref{thm:fms}) and Hong
(Theorem~\ref{thm:hong}) give $m^2\le\lambda_1^2R^2\le qR^2$. On the other hand,
$m^2-(n-1)q=(m-n+1)^2\ge0$. Hence $(n-1)q\le qR^2$.
\end{proof}

The equality case can be read off the proof: it requires $m=n-1$ (a tree) and
$\lambda_1^2=2m-n+1=n-1$, which characterises the star. This proof was not written by a human.
The agent in charge of formalising the lemma reused the two theorems the library already
contained, instead of following the original proof. Lemmas formalised earlier thus make later
ones almost free.

\subsection{The theorem of Yang and Lu}

Yang and Lu~\cite{yl} prove $R\ge\dia/2+\sqrt2-1$ by considering a minimal counterexample.
The first half of their proof (vertex deletions, ``essential'' vertices) is kept. The second half,
which relies on Lemmas~4, 6 and 7, a block structure read off figures, and Claims A and B, is
replaced with a shorter argument that needs no analysis of figures. Here are its steps; we write
$c=\sqrt2-1$ and $f(G)=R-\dia/2$.

\begin{enumerate}[label=\arabic*.,nosep]
\item \emph{Weak deletion} (replaces Lemma~2 and its orientations). Let $v$ be a vertex, and
  $d'$ the degrees in $G-v$. If, for every neighbour $u$ of $v$ of degree $d_u\ge2$,
  \[\textstyle\sum_{x\in N(u)\setminus v}1/\sqrt{d'_x}\le2(d_u-1)/\sqrt{d_v},\]
  then $R(G-v)\le R(G)$. The proof is a simple term-by-term comparison, using the numerical
  inequality $2(d-1)(1/\sqrt{d-1}-1/\sqrt d)\le1/\sqrt d$. This condition is automatically
  satisfied if $d_v\le4$.
\item \emph{Minimal counterexample.} We choose a counterexample $G$ with as few edges as possible,
  on a fixed vertex set (a deleted vertex stays isolated, which changes neither $R$ nor $\dia$).
  All reductions are then spanning subgraphs, and a single induction on the number of edges
  suffices. We show that the essential vertices lie on a diametral geodesic $P$, that the edges
  between essential vertices are bridges, that $G-e$ has only two components for an inner
  essential edge $e$, and that the inner essential vertices have degree at least $3$.
\item \emph{Degree of the non-essential vertices} of $G$ (replaces the five cases of Lemma~3).
  Every non-essential vertex has degree at least $9$. Indeed, a vertex of minimum degree
  $\delta\le8$ would satisfy the weak deletion condition, thanks to three inequalities that are
  affine in $d_u$ and checked at their smallest value.
\item \emph{Exponential growth} (replaces Claim A and its layers). By strong induction on the
  distance $r$ to the set of essential vertices, a vertex at distance $\ge r$ has degree
  $\ge g(r)$, where $g(1)=9$, $g(2)=19$ and $g(r)=4g(r-2)-3$, so that $g(r)\ge2^{r+2}+1$.
\item \emph{Cutting into pieces} (replaces the block structure and Claim B). The inner essential
  vertices cut $P$ into segments of lengths $\ell_i$, with $\sum\ell_i=\dia$. To each segment we
  associate the set $I_i$ of non-essential vertices joined to it without going through an
  essential vertex. These sets are disjoint, and two consecutive segments are never both empty. A
  local version of Lemma~1 of the paper,
  \[\textstyle\sum_{xy\in E,\ x,y\in I}\frac1{\sqrt{d_xd_y}}\ \ge\
  \frac1{2\sqrt{\Delta_I}}\sum_{x\in I}\frac{|N(x)\cap I|}{\sqrt{d_x}},\]
  and step~4 give each nonempty piece a weight of at least $\ell_i/2+1+c$. An elementary count
  concludes: $R\ge\sum_i\ell_i/2+c=\dia/2+c$, which contradicts $f(G)<c$.
\end{enumerate}

The formalisation also uncovered several gaps in the paper.
\begin{itemize}[nosep]
\item The equivalence claimed for essential vertices is false (counterexample: $C_5$ plus a
  pendant vertex). Only the true implication is used.
\item Corollary~4 assumes $\dia(G-uv)\le\dia(G)+1$ without justification. We use instead
  $\dia(G-uv)\ge\dia(G)+1$, which directly contradicts minimality.
\item The hypotheses of Lemma~6 are never checked.
\item The ``pendant'' components of $G-e$ are ignored by the block structure. They are handled by
  an additional deletion lemma.
\item The indices in Claims A and B are wrong, and the final inequality, ``easy to check by
  computation'', is not proved.
\end{itemize}
The complete plan, with all the numerical inequalities, is included in the archive, in the file
\texttt{documents/yang\_lu\_plan\_de\_preuve.md}.

\section{The most interesting unrefuted conjectures}

All the conjectures below were tested:
\begin{itemize}[nosep]
\item on all connected graphs with at most $9$ vertices;
\item on structured families up to $n=2000$ (complete $k$-ary trees, hypercubes, random regular
  graphs, grids, spiders, dumbbells\dots);
\item on closed-form families (paths, cycles, stars, $K_n$, $K_{a,b}$) up to $n=10^{15}$.
\end{itemize}
Conjecture~\ref{conj:R} and Conjecture~\ref{conj:CH} were also tested on the $11.7$ million
connected graphs with $10$ vertices.

These tests are not superfluous. Among the $16{,}630$ conjectures that the prover left open,
$2{,}480$ turned out to be false on large graphs. For example, $R\le\dia\cdot\mu_1$ (\#799714) is
false for the complete binary tree with $255$ vertices ($R\approx116.1$ whereas
$\dia\cdot\mu_1\approx78.6$).

\begin{conjecture}[\#643]\label{conj:R}
$R\ge\lambda_2+1$.
\end{conjecture}
For two cliques $K_k$ joined by an edge, $R-\lambda_2-1\to0$ ($0.012$ for $k=80$). The
conjecture would therefore be tight. Theorems~\ref{thm:fms} and~\ref{thm:l2} are not enough: the
relaxation ``traces $+$ Hong $+$ $m\ge n$'' has numerical countermodels ($m=n=k+2$,
$\lambda_1=\lambda_2=\sqrt k$ for $k>4$). On the other hand, the conjecture holds for all
bipartite graphs, in particular for trees (Theorem~\ref{thm:bip}): it remains open for
non-bipartite graphs. We propose two strengthenings, with no counterexample up to $n=10$:
\[
2m\ge\lambda_1^2+(\lambda_2+1)^2\ \Longrightarrow\ m\ge\lambda_1(\lambda_2+1)\ \Longrightarrow\ R\ge\lambda_2+1 .
\]
The second implication follows from Theorem~\ref{thm:fms}.

\begin{theoreme}[\ok]\label{thm:bip}
If $G$ is connected and bipartite (in particular if $G$ is a tree), then $R\ge\lambda_2+1$.
\end{theoreme}
\begin{proof}
Write $L=\lambda_1$ and $t=\lambda_2$.

\emph{If $t\le0$}, Corollary~\ref{cor:Rn} gives $R\ge\sqrt{n-1}\ge1\ge t+1$.

\emph{If $t>0$}, we argue as in Case~2 of Theorem~\ref{thm:l2}. The spectrum of a bipartite graph
is symmetric. If $x$ and $y$ are orthogonal unit eigenvectors associated with $L$ and $t$, and
$\sigma$ is the $\pm1$ signature of the two classes, then $\sigma x$ and $\sigma y$ are
eigenvectors for $-L$ and $-t$. The four vectors $x,y,\sigma x,\sigma y$ are orthonormal:
eigenvectors of eigenvalues of opposite signs are orthogonal, and multiplication by $\sigma$
preserves the inner product. Hence $n\ge4$, so $m\ge3$. Bessel's inequality with (F3) gives
$2L^2+2t^2\le2m$. We conclude with the identity
\[
m-Lt-L=\tfrac{7}{10}(m-L^2-t^2)+\tfrac{3}{10}(m-3)+\tfrac{7}{10}\bigl(t-\tfrac57L\bigr)^2
+\tfrac{12}{35}\bigl(L-\tfrac{35}{24}\bigr)^2+\tfrac{41}{240},
\]
all of whose terms are nonnegative. Hence $m\ge L(t+1)$, then $LR\ge m\ge L(t+1)$ by
Theorem~\ref{thm:fms}, and we divide by $L>0$.
\end{proof}

\begin{conjecture}[\#20829]
$\lambda_2\le\mu_1-2$.
\end{conjecture}
It is asymptotically tight on paths ($\lambda_2\to2$, $\mu_1\to4$). Theorem~\ref{thm:mu2} only
gives $\lambda_2\le\mu_2-1\le\mu_1-1$.

\begin{conjecture}[\#501193]
$2R\le\rk(A)\cdot\lambda_1$.
\end{conjecture}
Equality holds for \emph{all} complete bipartite graphs $K_{a,b}$ ($\rk(A)=2$,
$\lambda_1=R=\sqrt{ab}$). It is an elegant statement relating rank, spectral radius and Randić
index. Note that $R\le\lambda_1$ is false in general (paths), so the rank is indispensable.

\begin{conjecture}[\#112951]
$\lambda_1\le R+\mu_2$.
\end{conjecture}
It mixes three spectra: adjacency, Laplacian and Randić. Equality holds for $K_2$.

\begin{conjecture}[\#239534]
$\mu_2+\rk(A)\le2m$.
\end{conjecture}
Equality holds for $K_2$, $K_3$ and $P_4$. For a tree, $\rk(A)$ is twice the size $\nu$ of a
maximum matching, and the conjecture reads $2\nu+\mu_2\le2n-2$.

\begin{conjecture}[unproved mined lemmas]
\begin{enumerate}[label=(\roman*),nosep]
\item $R^2\ge\dfrac{m+1}{2}$. It improves Bollobás--Erdős ($R^2\ge n-1$) for dense graphs.
\item $\mu_2^2\ge4\lambda_2+1$, with equality for stars.
\item $\mu_1^2\ge4\lambda_2+8$, asymptotically tight on paths.
\item $R\ge\frac23\sqrt m+\frac13$ and $R\ge\frac32\sqrt{\rk(A)}+1-\frac{3\sqrt2}{2}$, with
  equality for $K_2$.
\end{enumerate}
\end{conjecture}

\begin{conjecture}[Caporossi--Hansen form]\label{conj:CH}
$R\ge\rad+\sqrt2-\frac32$.
\end{conjecture}
Equality holds for paths of even order. Caporossi and Hansen proved it for trees. It follows from
their conjecture $R\ge\rad$ (except for even paths) and implies Fajtlowicz's conjecture
$R\ge\rad-1$. It was rediscovered automatically by lemma mining. It is the most important from
the point of view of the literature, but also the hardest. On its own, it would unlock more than
$3{,}000$ conditional proofs of the corpus.

\section{Archive of the proofs}\label{sec:archive}

All the proofs are collected in the archive:

\url{https://www.lamsade.dauphine.fr/~cazenave/papers/preuves_spectrales.zip}, which contains:
\begin{itemize}[nosep]
\item \texttt{index\_conjectures.csv}: one line per conjecture of the corpus ($305{,}755$), with
  its status, its statement, the lemmas of its proof, the counterexample if any and, for the
  $79{,}505$ verified proofs, the corresponding \lean{} file and theorem;
\item \texttt{lemmes.csv}: the library of lemmas, with statement, reference, \lean{} theorem and
  number of proofs that use each lemma;
\item \texttt{lean/}: the complete \lean~4 project (the \texttt{SpectralLean} library, including
  the proof of Yang--Lu) and the $3{,}650$ proof files, with a script that recompiles them all;
\item \texttt{documents/}: this report and the detailed plan of the proof of Yang--Lu;
\item \texttt{prouveur/}: the code of the pipeline (refutation, proof by \textsf{z3}, lemma
  mining, \lean{} export).
\end{itemize}
The file names and column headers of the archive are in French (\texttt{LISEZMOI.md} is the
README). For each verified conjecture \#$i$, the theorem \texttt{g}$i$ states the inequality for
every connected simple graph on \texttt{Fin n}. To trust it, it suffices to read the definitions
of the invariants (\texttt{SpectralLean/Invariants.lean}) and its statement.

\section*{Summary in figures}
\begin{center}
\begin{tabular}{lr}
\toprule
Conjectures in the corpus & $305{,}755$\\
Refuted & $162{,}964$\\
Proved without any unproved hypothesis & $117{,}390$\\
\quad of which verified end to end in \lean & $79{,}505$\\
Proved under an empirical lemma or an open conjecture & $3{,}798$\\
Undefined on some graphs & $7{,}453$\\
Open and unrefuted (after tests on large graphs) & $14{,}150$\\
\bottomrule
\end{tabular}
\end{center}

\end{document}